# ON DIOPHANTINE EQUATION $X^2 = 2Y^4 - 1$

Florentin Smarandache
University of New Mexico
200 College Road
Gallup, NM 87301, USA
E-mail: smarand@unm.edu

Abstract: In this note we present a method of solving this Diophantine equation, method which is different from Ljunggren's, Mordell's, and R.K.Guy's.

In his book of unsolved problems Guy shows that the equation $x^2 = 2y^4 - 1$ has, in the set of positive integers, the only solutions $(1,1)$ and $(239,13)$; (Ljunggren has proved it in a complicated way). But Mordell gave an easier proof.

We'll note $t = y^2$. The general integer solution for $x^2 - 2t^2 + 1 = 0$ is

$$\begin{cases} x_{n+1} = 3x_n + 4t_n \\ t_{n+1} = 2x_n + 3t_n \end{cases}$$

for all $n \in \mathbb{N}$, where $(x_0, y_0) = (1, \varepsilon)$, with $\varepsilon = \pm 1$ (see [6]) or

$$\begin{pmatrix} x_n \\ t_n \end{pmatrix} = \begin{pmatrix} 3 & 4 \\ 2 & 3 \end{pmatrix}^n \cdot \begin{pmatrix} 1 \\ \varepsilon \end{pmatrix},$$ for all $n \in \mathbb{N}$, where a matrix to the power zero is equal to the unit matrix $I$.

Let's consider $A = \begin{pmatrix} 3 & 4 \\ 2 & 3 \end{pmatrix}$, and $\lambda \in \mathbb{R}$. Then $\det(A - \lambda \cdot I) = 0$ implies $\lambda_{1,2} = 3 \pm \sqrt{2}$, whence if $v$ is a vector of dimension two, then: $Av = \lambda_{1,2} \cdot v$.

Let's consider $P = \begin{pmatrix} 2 & 2 \\ \sqrt{2} & -\sqrt{2} \end{pmatrix}$ and $D = \begin{pmatrix} 3+2\sqrt{2} & 0 \\ 0 & 3-2\sqrt{2} \end{pmatrix}$. We have $P^{-1} \cdot A \cdot P = D$, or

$$A^n = P \cdot D^n \cdot P^{-1} = \begin{pmatrix} \frac{1}{2}(a+b) & \frac{\sqrt{2}}{2}(a-b) \\ \frac{\sqrt{2}}{4}(a-b) & \frac{1}{2}(a+b) \end{pmatrix},$$

where $a = \left(3 + 2\sqrt{2}\right)^n$ and $b = \left(3 - 2\sqrt{2}\right)^n$.

Hence, we find:

$$\begin{pmatrix} x_n \\ t_n \end{pmatrix} = \begin{pmatrix} \dfrac{1+\varepsilon\sqrt{2}}{2}\left(3+2\sqrt{2}\right)^n + \dfrac{1-\varepsilon\sqrt{2}}{2}\left(3-2\sqrt{2}\right)^n \\ \dfrac{2\varepsilon+\sqrt{2}}{4}\left(3+2\sqrt{2}\right)^n + \dfrac{2\varepsilon-\sqrt{2}}{4}\left(3-2\sqrt{2}\right)^n \end{pmatrix}, \ n \in \mathbb{N}.$$

Or $y_n^2 = \dfrac{2\varepsilon+\sqrt{2}}{4}\left(3+2\sqrt{2}\right)^n + \dfrac{2\varepsilon-\sqrt{2}}{4}\left(3-2\sqrt{2}\right)^n, \ n \in \mathbb{N}$.

For $n=0, \varepsilon=1$ we obtain $y_0^2 = 1$ (whence $x_0^2 = 1$), and for $n=3, \varepsilon=1$ we obtain $y_3^2 = 169$ (whence $x_3 = 239$).

(1) $\quad y_n^2 = \varepsilon \sum_{k=0}^{\left[\frac{n}{2}\right]} \binom{n}{2k} \cdot 3^{n-2k} 2^{3k} + \sum_{k=0}^{\left[\frac{n-1}{2}\right]} \binom{n}{2k+1} \cdot 3^{n-2k-1} 2^{3k+1}$

We still must prove that $y_n^2$ is a perfect square if and only if $n = 0, 3$.

We can use a similar method for the Diophantine equation $x^2 = Dy^4 \pm 1$, or more generally: $C \cdot X^{2a} = DY^{2b} + E$, with $a, b \in \mathbb{N}^*$ and $C, D, E \in \mathbb{Z}^*$; denoting $X^a = U$, $Y^b = V$, and applying the results from F.S. [6], the relation (1) becomes very complicated.